# Every rearrangement-invariant quasi-Banach function space is an interpolation space between two Lorentz spaces

Leo R. Ya. Doktorski[1]

## Abstract.

We extend some results of Kwok-Pun Ho. In particular, it will be shown that every rearrangement-invariant quasi-Banach function space $E$ on a totally σ-finite measure space Ω with a non-atomic measure $\mu$ can be expressed in the form

$$E = \mathcal{F}\left(L_{p_0,q_0}, L_{p_1,q_1}\right)$$

for an interpolation functor $\mathcal{F}$, where the construction of the functor $\mathcal{F}$ is given based on the space $E$ itself.

Keywords: Rearrangement-invariant quasi-Banach function spaces, Boyd indices, Lorentz spaces.

Mathematics Subject Classification: 46A16 ▪ 46B70 ▪ 46E30

## 1. Introduction

Let $E$ be a rearrangement-invariant quasi-Banach function space on $\mathbb{R}^n$ with the Boyd indices $p_E$ and $q_E$. In [12, Corollary 4.3 and 15, Corollary 3.4], Kwok-Pun Ho has shown that under the conditions

$$0 < p_0 < p_E \leq q_E < p_1 < \infty$$

the space $E$ can be expressed in the form

$$E = \mathcal{F}\left(L_{p_0}, L_{p_1}\right)$$

for an interpolation functor $\mathcal{F}$. Moreover, the construction of the functor $\mathcal{F}$ is given based on the space $E$ itself and Peetre's $K$-functional $K\left(t, f; L_{p_0}, L_{p_1}\right)$.

The main purpose of this short note is to extend this result. We consider rearrangement-invariant quasi-Banach function spaces on a totally $\sigma$-finite measure space Ω with a non-atomic measure $\mu$ such that $\mu(\Omega) = L$ where $L \in \{1, \infty\}$. It will be shown that under the conditions

$$0 < q_0, q_1 \leq \infty \text{ and } 0 < p_0 < p_E \leq q_E < p_1 \leq \infty \text{ or } 0 < p_0 < p_E \leq q_E = p_1 = q_1 = \infty$$

each such space $E$ can be expressed in the form

$$E = \mathcal{F}\left(L_{p_0,q_0}, L_{p_1,q_1}\right)$$

for an interpolation functor $\mathcal{F}$ introduced by Kwok-Pun Ho. Here $L_{p_i,q_i}$ $(i \in \{0,1\})$ are the Lorentz spaces. Note that the case $q_E = \infty$ is also allowed. This leads to an important corollary: for every rearrangement-invariant Banach function space $E$ there exist $p_0$ and $p_1$ such that $0 < p_0 < p_1 \leq \infty$ and

$$E = \mathcal{F}\left(L_{p_0}, L_{p_1}\right).$$

This paper is organized as follows. Section 2 contains some notations, definitions and technical results. In Section 3, we will recall some basic definitions and constructions related to interpolation theory and $K$-method of real interpolation. In Subsection 3.3, we will recall the definition of the interpolation functor $\mathcal{F}$ introduced by Kwok-Pun Ho and prove an analogue of Theorem 4.1 from [12] which includes the case $q_E = \infty$. The main result of this paper is proved in Section 4. Some important corollaries of it will be formulated and compared with known results in Section 5.

[1] Fraunhofer IOSB (Institute of Optronics, System Technologies and Image Exploitation), Department Object Recognition, Gutleuthausstr. 1, 76275 Ettlingen, Germany. ORCID iD: 0000-0001-5977-1000
Email: doktorskileo@gmail.com; leo.doktorski@iosb.fraunhofer.de

## 2. Preliminaries

For $f$ and $g$ being positive functions, we write $f \preccurlyeq g$ if $f \leq Cg$, where the positive constant $C$ is independent of all significant quantities. Two positive functions $f$ and $g$ are considered equivalent ($f \sim g$) if $f \preccurlyeq g$ and $g \preccurlyeq f$.

We adopt the conventions $1/\infty = 0$ and $1/0 = \infty$.

A few more denotations:

$\|*\|_{p,(a,b)}$ is the usual (quasi-)norm on the Lebesgue space $L_p \equiv L_p(a,b)$ on the interval $(a,b)$ ($0 < p \leq \infty$, $0 \leq a < b \leq \infty$). If $a = 0$ and $b = \infty$, we write $\|*\|_p$.

$f_n \nearrow f$: the sequence $\{f_n\}$ $(n = 1,2,\ldots)$ of real-valued non-negative functions is monotonically increasing and we have the pointwise convergence $f_n \to f$ a.e..

Let $\omega \subset (0,\infty)$. We use $\chi_\omega$ to denote the characteristic function of $\omega$.

### 2.1 Quasi-normed spaces

We recall some basic facts concerning quasi-normed spaces. As a reference for this subject, we recommend [4].

Definition 1. A quasi-norm $\|x\|_X$ defined on a linear vector space $X$ (over a field $\mathbb{K}$ of real or complex numbers) is a map $X \to \mathbb{R}^+$ such that

(i) $\|x\|_X = 0 \iff x = 0$,

(ii) $\|\lambda x\|_X = |\lambda| \|x\|_X$ for all $\lambda \in \mathbb{K}$ and $x \in X$,

(iii) $\|x + y\|_X \leq C_X(\|x\|_X + \|y\|_X)$ for all $x, y \in X$, where the constant $C_X$ is independent of $x, y$.

A space $X$ with a quasi-norm $\|x\|_X$ is called a quasi-normed space. When it is complete, it will be called quasi-Banach space. The inequality (iii) above is called the "*quasi-triangle inequality*" and the constant $C_X$ is called "*quasi-triangle constant*". Obviously, $C_X \geq 1$. If $C_X = 1$, then $X$ is a normed space.

Let $X$ and $Y$ be two quasi-normed spaces with quasi-norms $\|*\|_X$ and $\|*\|_Y$, respectively. Throughout the paper, we write $X \hookrightarrow Y$ to indicate that $X$ is continuously embedded in $Y$, i.e. there exists a constant $C > 0$ such that $\|x\|_Y \leq C\|x\|_X$ for every $x$ in $X$. We write $X = Y$ if $X \hookrightarrow Y$ and $Y \hookrightarrow X$. In this case, $X$ and $Y$ are equal as sets and as linear spaces and their quasi-norms are equivalent ($\|x\|_X \sim \|x\|_Y$). We write $T: X \to Y$ to indicate that $T$ is a bounded operator from $X$ to $Y$. Its norm we denote by $\|T\|_{X \to Y}$.

T. Aoki [1] and S. Rolewicz [33] characterized quasi-Banach spaces as follows. See, e.g., [26] and [4, Lemmas 3.10.1, 3.10.2].

Theorem 2. (Aoki-Rolewicz Theorem). Let $X$ be a quasi-normed space with the quasi-norm $\|x\|_X$ and the quasi-triangle constant $C_X$. Then for $\kappa = \frac{1}{\log_2(2C_X)}$ there exists an equivalent quasi-norm $\|x\|_*$ on $X$ such that

$$\|x\|_* \leq \|x\|_X \leq 2C\|x\|_* \qquad \text{and} \qquad \|x + y\|_*^\kappa \leq \|x\|_*^\kappa + \|y\|_*^\kappa.$$

## 2.2 Rearrangement-invariant quasi-Banach function spaces

Banach function spaces were introduced in 1955 by Luxemburg. They include not only Lebesgue and Lorentz spaces mentioned above but also Zygmund, Lorentz–Zygmund and generalised Lorentz–Zygmund spaces. Various definitions of Banach and quasi-Banach function spaces exist in the literature. See, for example, [6, 26, 30], [3, Chapter 1, Definition 1.1], [8, Definition 3.1] and [29, Definition 2.10.], and references therein. A discussion of the various definitions of a (quasi-) Banach functional space can be found in [24]. A little different approach via so-called ideal spaces can be found, e.g., in [19, 20, 22].

Throughout the whole paper, $L \in \{1, \infty\}$. Let $(\Omega, \mu)$ be a totally $\sigma$-finite measure space with a non-atomic measure $\mu$ such that $\mu(\Omega) = L$. Note that this measure space is resonant. (See [3, Chapter 2, Definition 2.3, p.45 and Theorem 2.7, p 52].) By $\mathfrak{M}(\Omega, \mu)$ we will denote the set of all complex-valued, $\mu$ −measurable and finite almost everywhere (abbreviation: a.e.) functions on $\Omega$. As is customary, we will identify functions that coincide a.e.. By $\mathfrak{M} \equiv \mathfrak{M}(0, L)$ we will denote the set of all complex-valued, Lebesgue measurable and finite almost everywhere functions on $(0, L)$.

The *non-increasing rearrangement* $f^*$ of $f \in \mathfrak{M}(\Omega, \mu)$ is defined by (see, e.g., [3])

$$f^*(t) = \{\lambda > 0: \mu\{x \in \Omega: |f(x)| > \lambda\} \leq t\}, \qquad t > 0.$$

It is known that the non-increasing rearrangement is unique. The functions $f$ and $g$ are called equimeasurable if $f^* = g^*$. More information about the non-increasing rearrangement $f^*$ can be found, e.g., in [3, Chapter 2]. We recall the definition of rearrangement-invariant quasi-Banach function space [12, 13, 28, 30].

Definition 3. A quasi-Banach space $E \subset \mathfrak{M}(\Omega, \mu)$ is called a *rearrangement-invariant* (abbreviation *r.i.*) *quasi-Banach function space* if there exists a mapping $\rho_E: \mathfrak{M}(0, L) \to [0, \infty]$ satisfying

(P1) $\rho_E$ is a quasi-norm,

(P2) $\rho_E$ has the *lattice property*: if $f, g \in \mathfrak{M}(\Omega, L)$ and $|g| \leq |f|$ a. e., then $\rho_E(g) \leq \rho_E(f)$,

(P3) $\rho_E$ has the *Fatou property*: if $0 \leq f_n \nearrow f$ a.e. $(n = 1,2, \dots)$, then $\rho_E(f_n) \nearrow \rho_E(f)$,

(P4) for all measurable sets $\omega \subset (0, L)$ with finite measure $\rho_E(\chi_\omega) < \infty$,

so that

$$\|f\|_E = \rho_E(f^*), \qquad \forall f \in E.$$

It is known that all r.i quasi-Banach function spaces are complete. See [8, Lemma 3.6], [30, Corollary 3.8] and [24, Proposition 2.2]. Therefore, the use of the term "Banach" in the name of these spaces is justified.

Definition 3 includes the Fatou property (P3). Many authors (see, for example, [9, 10]) have investigated a broader class of the so-called *symmetric (quasi-) Banach function spaces*. The definition of these spaces does not include the Fatou property. The classical monographs [21] and [3] are devoted to symmetric and to r.i. Banach function spaces, respectively. We limit ourselves to considering only r.i. (quasi-) Banach function spaces.

For simplicity, throughout this paper we consider functions from $\mathfrak{M}(0, L)$. Because of the Luxemburg representation theorem (See Theorem 3.1 and Proposition 3.3 in [29].) all our results can be reformulated for $\mu$-measurable functions defined on a $\sigma$-finite measure space $(\Omega, \mu)$ with a non-atomic measure $\mu$.

## 2.3 Lorentz spaces

Lorentz spaces, introduced by George G. Lorentz [23], are generalisations of the more familiar Lebesgue space $L_p$. See, e.g., [17].

Definition 4. Let $0 < p, q \leq \infty$. The Lorentz spaces $L_{p,q} \equiv L_{p,q}(0, L)$ are the sets of all $f \in \mathfrak{M}(0, L)$ with the finite (quasi-)norms

$$\|f\|_{p,q} \coloneqq \left\| t^{\frac{1}{p}-\frac{1}{q}} f^*(t) \right\|_{q,(0,L)} \text{ if } 0 < q < \infty,$$

$$\|f\|_{p,\infty} \coloneqq \sup_{0<t<L} t^{\frac{1}{p}} f^*(t) \qquad \text{if } q = \infty.$$

The following lemma establishes the basic properties of Lorentz spaces that we need.

Lemma 5. Let $0 < p, q \leq \infty$. The Lorentz space $L_{p,q}$ is a r.i. quasi-Banach function space, if and only if one of the following conditions is satisfied:

(i) $p < \infty$;

(ii) $p = q = \infty$.

Furthermore, if these conditions are not met, then $L_{p,q} = \{0\}$.

## 2.4 Dilation operators and Boyd indices

As a reminder, for simplicity, we will henceforth consider functions from $\mathfrak{M}(0, L)$.

An important property of r.i. quasi-Banach function spaces is that the dilation operator is bounded on these spaces. This makes it possible to define the Boyd indices. They were introduced by Boyd in [5 and 6] and are useful in various problems concerning rearrangement-invariant spaces. To readers interested in this topic we would recommend [26] and monographs [3, Chapter 3, Section 5] and [21, Chapter 2]. New results can be found in [9, 16, 29, 30, 32] and references therein. See also [10, 20, 22, 25].

Definition 6. Let $0 < a < \infty$. The dilation operator on $\mathfrak{M}(0, L)$ we define by

$$(D_a f)(t) := f(at)\chi_{(0,L)}(at), \qquad 0 < t < L.$$

Obviously,

$$D_a(\lambda f) = \lambda D_a(f), D_a(f + g) = D_a(f) + D_a(g) \text{ and } D_a f^* = (D_a f)^*.$$

Let $E$ be some r.i. quasi-Banach function space. We are interested in the behavior of the norm of the dilation operator in $E$: $\|D_a\| \equiv \|D_a\|_{E\to E}$. Because $f^*(at) \leq f^*(t)$ for $a \geq 1$, $\|D_a\| \leq 1$. It is known (see [30, Lemma 3.18 and Theorem 3.23], [9, Lemma 4.4], [7, 21]) that the dilation operators $D_a$ are bounded for every choice of $a \in (0, \infty)$. Additionally, the function $a \to \|D_a\|$ is submultiplicative. This means that

$$\|D_{ab}\| \leq \|D_a\| \|D_b\| \qquad \text{for all } 0 < a, b < \infty.$$

All the main results of this paper are formulated in terms of Boyd indices. We use the definition of the lower and upper Boyd indices in the following form. See, e.g., [9, 10, 13, 15, 28].

Definition 7. Let $E$ be a r.i. quasi-Banach function space. Define the lower Boyd index $p_E$ and the upper Boyd index $q_E$ of $E$ by

$$p_E = \sup\left\{p > 0 : \exists c > 0 \text{ such that } \forall 0 < a < 1 \ \|D_a\| \leq c a^{-\frac{1}{p}}\right\},$$

$$q_E = \inf\left\{q > 0 : \exists c > 0 \text{ such that } \forall a \geq 1 \ \|D_a\| \leq c a^{-\frac{1}{q}}\right\},$$

respectively.

It is known (see, e.g., [26], [9, Corollary 4.7]) that

$$0 < p_E \leq q_E \leq \infty$$

and if $E$ is a r.i. Banach function space, then

$$1 \leq p_E \leq q_E \leq \infty.$$

In [9, Proposition 4.6] it is shown that

$$p_E = \sup_{s>1} \frac{\ln s}{\ln\|D_{1/s}\|} = \lim_{s\to\infty} \frac{\ln s}{\ln\|D_{1/s}\|} \qquad \text{and} \qquad q_E = \sup_{0<s<1} \frac{\ln s}{\ln\|D_{1/s}\|} = \lim_{s\to 0+} \frac{\ln s}{\ln\|D_{1/s}\|}.$$

In many texts (see, e.g., [3, Chapter 3, Definitions 5.10 and 5.12]), the lower and upper Boyd indices of $E$ are alternatively defined as the quantities

$$\underline{\alpha}_E = \sup_{0<s<1} \frac{\ln\|D_{1/s}\|}{\ln s}, \qquad \text{and} \qquad \overline{\alpha}_E = \sup_{s>1} \frac{\ln\|D_{1/s}\|}{\ln s}.$$

It is clear that $\underline{\alpha}_E = \frac{1}{q_E}$ and $\overline{\alpha}_E = \frac{1}{p_E}$.

### 2.5 Averaging operators

In this subsection, we define some averaging operators and formulate their properties which we need. They play an important role below. We are interested in the conditions under which these operators are bounded in r.i. quasi-Banach function spaces. The main result of this subsection is Theorem 11. It is a corollary of Theorem 2 from [28].

Definition 8. Let $0 < U, V, W < \infty$. We consider the following averaging (Hardy-type) operators on $\mathfrak{M}(0, L)$.

$$\left(H^{(U,W)} f\right)(t) := t^{-\frac{1}{U}} \left\{ \int_0^t \left[ v^{\frac{1}{U}} f^*(v) \right]^W \frac{dv}{v} \right\}^{\frac{1}{W}},$$

$$\left(H^{(U,\infty)} f\right)(t) := t^{-\frac{1}{U}} \sup_{0<v<t} v^{\frac{1}{U}} f^*(v),$$

$$\left(H_{(V,W)} f\right)(t) := t^{-\frac{1}{V}} \left\{ \int_t^L \left[ v^{\frac{1}{V}} f^*(v) \right]^W \frac{dv}{v} \right\}^{\frac{1}{W}},$$

and

$$\left(H_{(V,\infty)} f\right)(t) := t^{-\frac{1}{V}} \sup_{t<v<L} v^{\frac{1}{V}} f^*(v).$$

Without stipulating this each time, we will assume that

$$0 < U, V < \infty \qquad \text{and} \qquad 0 < W \le \infty.$$

If necessary, we will impose additional constraints on these parameters. These operators are studied in [22, 28]. In [5], in [6, (19) and (20)], and in [3, Definition III.5.14], averaging operators

$$(P_\alpha f)(t) = t^{-\alpha} \int_0^t v^\alpha f^*(v) \frac{dv}{v}, \qquad (0 < \alpha \le 1),$$

and

$$(Q_\alpha f)(t) = t^{-\alpha} \int_t^\infty v^\alpha f^*(v) \frac{dv}{v}, \qquad (0 \le \alpha < 1),$$

are considered. Obviously, $P_\alpha = H^{(\alpha^{-1},1)}$ and $Q_\alpha = H_{(\alpha^{-1},1)}$. Slightly different averaging operators are studied in [25]. For $0 < U < \infty$ put

$$f^{**}_{(U)}(t) := \left( t^{-1} \int_0^t f^*(v)^U dv \right)^{\frac{1}{U}}.$$

These averaging operators were used in [34 and 36]. Cf. also [11, (4.2)]. Obviously,

$$f^{**}_{(U)}(t) = \left(H^{(U,U)} f\right)(t).$$

In particular,

$$f^{**}(t) := t^{-1} \int_0^t f^*(v) dv = f^{**}_{(1)}(t) = \left(H^{(1,1)} f\right)(t).$$

The operators $H^{(U,U)}$ and $H_{(V,V)}$ were also used in [7].

The next theorem plays a crucial role in our research below. Therein conditions are given for the operators $H^{(U,W)}$ and $H_{(V,W)}$ to be bounded in a r.i. quasi-Banach function space in terms of its Boyd indices. For the case $L = \infty$ it coincides with [28, Theorem 2] (Cf. [3, Theorem III.5.17].). Theorem 17 from [22] extends it for the case $L = 1$. However, there are no statements (ii) and (iv) in [22], which can also be easily proved.

Theorem 9. If $E$ is a r.i. quasi-Banach function space on $(0, L)$ then we have the following.

(i) If $0 < W < \infty$, then the operator $H^{(U,W)}$ is bounded from $E$ to $E$ if and only if $p_E > U$.

(ii) The operator $H^{(U,\infty)}$ is bounded from $E$ to $E$ if $p_E > U$.

(iii) If $0 < W < \infty$, then the operator $H_{(V,W)}$ is bounded from $E$ to $E$ if and only if $q_E < V$.

(iv) The operator $H_{(V,\infty)}$ is bounded from $E$ to $E$ if $q_E < V$.

Note that the reverse implications are not true in points (ii) and (iv) (see [28]).

The following lemma describes some properties of the averaging operators from Definition 8.

Lemma 10. Consider the averaging operators $H^{(U,W)}$ or $H_{(V,W)}$ from Definition 8. Then, for all $f \in \mathfrak{M}(0, L)$ and $t > 0$, the following holds:

(i) $\left(H^{(U,W)}f\right)(t) \succcurlyeq \left(H^{(U,\infty)}f\right)(t) \ge f^*(t)$;

(ii) $\left(H_{(V,W)}f\right)(t) \succcurlyeq f^*(2t)$ and $\left(H_{(V,\infty)}f\right)(t) \ge f^*(t)$.

Proof. First, we will prove assertion (i). Because $f^*$ is non-increasing, using the fact that $\left\{\int_0^v s^{\frac{W}{U}} \frac{ds}{s}\right\}^{\frac{1}{W}} \sim v^{\frac{1}{U}}$ for all $v > 0$, we have

$$f^*(t) = f^*(t)t^{-\frac{1}{U}} \sup_{0<v<t} v^{\frac{1}{U}} \le t^{-\frac{1}{U}} \sup_{0<v<t} v^{\frac{1}{U}} f^*(v) = \left(H^{(U,\infty)}f\right)(t)$$

$$\sim t^{-\frac{1}{U}} \sup_{0<v<t} f^*(v) \left\{\int_0^v s^{\frac{W}{U}} \frac{ds}{s}\right\}^{\frac{1}{W}} \le t^{-\frac{1}{U}} \sup_{0<v<t} \left\{\int_0^v \left(s^{\frac{1}{U}} f^*(s)\right)^W \frac{ds}{s}\right\}^{\frac{1}{W}}$$

$$= t^{-\frac{1}{U}} \left\{\int_0^t \left(s^{\frac{1}{U}} f^*(s)\right)^W \frac{ds}{s}\right\}^{\frac{1}{W}} = \left(H^{(U,W)}f\right)(t).$$

Let us prove (ii). Using that $f^*$ is non-increasing, we get

$$\left(H_{(V,W)}f\right)(t) \ge t^{-\frac{1}{V}} \left\{\int_t^{2t} \left(v^{\frac{1}{V}} f^*(v)\right)^W \frac{dv}{v}\right\}^{1/W}$$

$$\ge t^{-\frac{1}{V}} f^*(2t) \left\{\int_t^{2t} v^{\frac{W}{V}} \frac{dv}{v}\right\}^{1/W} \sim t^{-\frac{1}{V}} f^*(2t) t^{\frac{1}{V}} = f^*(2t).$$

The last inequality can be proved in a similar way. ∎

Lemma 10 allows one to prove the main result in this subsection.

Theorem 11. Let $E$ be a r.i. quasi-Banach function space with the quasi-norm $\|f^*\|_E$.

(i) If $0 < U < p_E$, then for all $f \in \mathfrak{M}$

$$\|f\|_E \sim \left\|H^{(U,W)} f\right\|_E.$$

(ii) If $q_E < V < \infty$, then for all $f \in \mathfrak{M}$

$$\|f\|_E \sim \left\|H_{(V,W)}\right\|_E.$$

Proof. First, note that if $H$ is one of the averaging operators $H^{(U,W)}$ from Definition 8 then by Lemma 10 (i), for all $f \in \mathfrak{M}$ and $t > 0$

$$f^*(t) \preccurlyeq (Hf)(t).$$

If $\|Hf\|_E < \infty$, the lattice property (P2) implies now

$$\|f\|_E \preccurlyeq \|Hf\|_E. \tag{1}$$

Similarly, if $H$ is one of the averaging operators $H_{(V,W)}$ from Definition 8 and $\|Hf\|_E < \infty$, then by Lemma 10 (ii), we observe that

$$\|f^*(2t)\|_E \preccurlyeq \|(Hf)(t)\|_E.$$

Because the operator $D_{\frac{1}{2}}$ is bounded in $E$, we have

$$\|f^*(t)\|_E = \left\|\left(D_{\frac{1}{2}}f^*\right)(2t)\right\|_E \preccurlyeq \|f^*(2t)\|_E.$$

In this way, one also gets assessment (1). Thus, one only needs to prove the inverse inequalities under the assumption $\|f\|_E < \infty$.

We continue the proof of (i). Part (ii) can be considered analogously. Let $H$ be one of the averaging operators $H^{(U,W)}$ from Definition 8. Because $0 < U < p_E$, by Theorem 9 (i) and (ii), we have

$$\|Hf\|_E \preccurlyeq \|f\|_E.$$

Thus, considering (1),

$$\|f\|_E \sim \|Hf\|_E.$$

■

## 3. Fundamentals of interpolation theory

In this section, we will recall some basic definitions and constructions related to the interpolation theory. For more details we refer to [3, 4, 35].

### 3.1 Basic definitions

In the following, $\bar{X} \equiv (X_0, X_1)$ will be a *compatible couple* of quasi-normed spaces. This means that $X_0$ and $X_1$ are continuously embedded in some common Hausdorff topological vector space. As usual,

$$X_0 + X_1 = \{\, f = f_0 + f_1 : f_0 \in X_0, f_1 \in X_1 \,\}.$$

The intersection $X_0 \cap X_1$ and the sum $X_0 + X_1$ are also quasi-norms spaces with the quasi-norms

$$\|f\|_{X_0 \cap X_1} = \max\{\|f\|_{X_0}, \|f\|_{X_1}\}$$

and

$$\|f\|_{X_0 + X} = \inf\{\|f_0\|_{X_0} + \|f_1\|_{X_1}; f_0 \in X_0, f_1 \in X_1 : f = f_0 + f_1\},$$

respectively.

A quasi-Banach space $X$ is called an *intermediate space* for the couple $\bar{X}$ (or an *intermediate space between* $X_0$ and $X_1$) if

$$X_0 \cap X_1 \hookrightarrow X \hookrightarrow X_0 + X_1.$$

Consider a linear operator $T: X_0 + X_1 \to X_0 + X_1$ whose restrictions to $X_0$ and $X_1$ are bounded from $X_i$ to $X_i$, $(i = 0,1)$. The space $X$ is called an *interpolation space* between $X_0$ and $X_1$ (or *with respect to* $\bar{X}$) if

$$T: X_0 + X_1 \to X_0 + X_1 \qquad \text{implies} \qquad T: X \to X$$

for each such linear operator $T$. If in addition,

$$\|T\|_{X \to X} \leq \max\{\|T\|_{X_0 \to X_0}, \|T\|_{X_1 \to X_1}\},$$

one says that $X$ is an *exact interpolation space*.

### 3.2 $K$-Method of real interpolation

We will use one of the most important ways of constructing interpolation spaces based on the use of *Peetre's K-functional*. For an arbitrary compatible couple $\bar{X} = (X_0, X_1)$ of quasi-normed spaces and for each $t > 0$, Peetre's $K$-functional is defined on $X_0 + X_1$ and given by

$$K(t, x; \bar{X}) \equiv K(t, x) := \inf\left(\|x_0\|_{X_0} + t\|x_1\|_{X_1} : x = x_0 + x_1, x_i \in X_i, i = 0,1\right).$$

For a fixed $x \in X_0 + X_1$ the function $t \to K(t, x)$ is continuous, non-decreasing, concave and non-negative on $(0, \infty)$, the function $t \to t^{-1}K(t, x)$ is continuous and non-increasing. Additional properties of the $K$-functional can be found, for example, in [3, 4, 35]. Below we assume that they are known to the reader. The following lemma can be proved in the usual way. See [31, Proposition 3.2].

Lemma 12. If $(X_0, X_1)$ is a compatible couple of quasi-normed spaces with the quasi-triangle constants $C_{X_0}$ and $C_{X_1}$ then

$$K(t, x_0 + x_1) \leq \max\{C_{X_0}, C_{X_1}\}\{K(t, x_0) + K(t, x_1)\}.$$

### 3.3 The interpolation functor associated with a r.i. quasi-Banach function space

We use the interpolation functor introduced by Kwok-Pun Ho. See [12, Definition 4.2], [15, Definition 3.1] and [13, Definition 2.1]. As shown in [12, 13, 14, 15, 16], this interpolation functor gives a number of applications on the boundedness of sublinear operators on rearrangement-invariant quasi-Banach function spaces.

Definition 13. Let $0 < \Theta, r < \infty$ and $E$ be a r.i. quasi-Banach function space on $(0, L)$ with the quasi-norm $\rho_E$. Let $(X_0, X_1)$ be a compatible couple of (quasi-) normed spaces. The space $\bar{X}_{\Theta,r;E} \equiv (X_0, X_1)_{\Theta,r;E}$ consists of all $x \in X_0 + X_1$ such that

$$\|x\|_{\bar{X}_{\Theta,r;E}} := \rho_E\left(t^{-\frac{1}{r}}K\left(t^{\frac{1}{\Theta}}, x, X_0, X_1\right)\right) < \infty.$$

Traditionally, the main parameter of a real interpolation method is denoted by the lowercase letter $\theta$ and satisfies the restrictions $0 < \theta < \infty$ (or $0 \leq \theta \leq \infty$). In the functor introduced by Ho, the main parameter lies between 0 and $\infty$. We use the capital letter $\Theta$ to indicate that the functor introduced by Ho is being used.

In the following theorem, compared to Theorem 4.1 from [12], the parameter $q_E$ can take the value $\infty$. In addition, another estimate for operators $T: (X_0, X_1) \to (Y_0, Y_1)_{\Theta,r;E}$ has been proven.

Theorem 14. Let $E$ be a r.i. quasi-Banach function space on $(0, L)$ with the quasi-norm $\rho_E$, $0 < \Theta < \infty$ and $0 < r < p_E$. Additionally, assume that

$$\frac{1}{\Theta} + \frac{1}{q_E} > \frac{1}{r} \qquad \text{if } q_E < \infty$$

or

$$\frac{1}{\Theta} \geq \frac{1}{r} \qquad \text{if } q_E = \infty.$$

Then $(*,*)_{\Theta,r;E}$ is an exact interpolation functor and if $(X_0, X_1)$ is a compatible couple of quasi-normed spaces, then the space $(X_0, X_1)_{\Theta,r;E}$ is also a quasi-normed space. Moreover, if the spaces $X_0, X_1$ and $E$ are normed, then $\bar{X}_{\Theta,r;X}$ is also a normed space.

Proof. (Cf. with the proof of Theorem 4.1 from [12].) Let $\bar{X} = (X_0, X_1)$ be a compatible couple of quasi-normed spaces. Denote the quasi-triangle constants of the spaces $E$, $X_0$ and $X_1$ by $C_E$, $C_{X_0}$ and $C_{X_1}$, respectively.

Step 1. To prove that $\bar{X}_{\Theta,r;E}$ is a quasi-normed linear space, it suffices to check the triangle inequality. If $x, y \in \bar{X}_{\Theta,r;E}$, by Lemma 12, we have

$$\|x+y\|_{\bar{X}_{\Theta,r;E}} = \rho_E\left(t^{-\frac{1}{r}}K\left(t^{\frac{1}{\Theta}}, x+y\right)\right)$$

$$\leq \rho_E\left(t^{-\frac{1}{r}}\max\{C_{X_0}, C_{X_1}\}\left\{K\left(t^{\frac{1}{\Theta}}, x\right)+K\left(t^{\frac{1}{\Theta}}, y\right)\right\}\right)$$

$$\leq \max\{C_{X_0}, C_{X_1}\}\, C_E\left\{\rho_E\left(t^{-\frac{1}{r}}K\left(t^{\frac{1}{\Theta}}, x\right)\right)+\rho_X\left(t^{-\frac{1}{r}}K\left(t^{\frac{1}{\Theta}}, y\right)\right)\right\}$$

$$= \max\{C_{X_0}, C_{X_1}\}\, C_E\{\|x\|_{\bar{X}_{\Theta,r;E}}+\|y\|_{\bar{X}_{\Theta,r;E}}\}.$$

The last estimate implies that if the spaces $X_0, X_1$ and $E$ are normed ($C_{X_0} = C_{X_1} = C_E = 1$), then $\bar{X}_{\Theta,r;E}$ is also a normed space.

Step 2. Let us show that $\bar{X}_{\Theta,r;E}$ is an intermediate space for the couple $\bar{X}$. The embedding $\bar{X}_{\Theta,r;E} \hookrightarrow X_0 + X_1$ is obvious because

$$\min\{1,t\}\,\|x\|_{X_0+X_1} \leq K(t,x).$$

Indeed, we have

$$\|x\|_{\bar{X}_{\Theta,r;E}} = \rho_E\left(t^{-\frac{1}{r}}K\left(t^{\frac{1}{\Theta}}, x\right)\right) \geq \rho_E\left(t^{-\frac{1}{r}}\min\left\{1, t^{\frac{1}{\Theta}}\right\}\right)\|x\|_{X_0+X_1}.$$

Thus, if $0 < \|x\|_{\bar{X}_{\Theta,r;E}} < \infty$, then $0 < \rho_E\left(t^{-\frac{1}{r}}\min\left\{1, t^{\frac{1}{\Theta}}\right\}\right) < \infty$ and

$$\|x\|_{X_0+X_1} \leq \left\{\rho_E\left(t^{-\frac{1}{r}}\min\left\{1, t^{\frac{1}{\Theta}}\right\}\right)\right\}^{-1}\|x\|_{\bar{X}_{\Theta,r;E}}.$$

This implies $\bar{X}_{\Theta,r;X} \hookrightarrow X_0 + X_1$. Let us prove that $X_0 \cap X_1 \hookrightarrow \bar{X}_{\Theta,r;X}$. By the definition of the $K$-functional, for any $x \in X_0 \cap X_1$, it holds

$$K(t,x) \leq \min\{t,1\}\,\|x\|_{X_0\cap X_1}.$$

So,

$$\|\mathrm{x}\|_{\bar{X}_{\Theta,r;E}} \leq \rho_E\left(t^{-\frac{1}{r}}\min\left\{t^{\frac{1}{\Theta}}, t\right\}\right)\|x\|_{X_0\cap X_1}.$$

Thus, it is enough to show that

$$\rho_E\left(\min\left\{t^{\frac{1}{\Theta}-\frac{1}{r}}, t^{-\frac{1}{r}}\right\}\right) < \infty.$$

Since $\rho_E$ is a quasi-norm, Aoki–Rolewicz Theorem 2 assures that there exists an equivalent quasi-norm (which we will also denote as $\rho_E$) and a positive number $\kappa_E$ such that $\rho_E^{\kappa_E}$ satisfies the triangle inequality. Hence, because $\frac{1}{\Theta} > 0$, we have

$$\rho_E^{\kappa_E}\left(\min\left\{t^{-\frac{1}{r}}, t^{\frac{1}{\Theta}-\frac{1}{r}}\right\}\right) \leq \rho_E^{\kappa_E}\left(t^{\frac{1}{\Theta}-\frac{1}{r}}\chi_{(0,1]}(t)\right)+\rho_E^{\kappa_E}\left(t^{-\frac{1}{r}}\chi_{(1,\infty)}(t)\right) := I + II.$$

In [12, proof of Theorem 4.1] it is shown that $II < \infty$. The proof uses only the condition $r < p_E$ and allows the value $p_E = \infty$. Note that if $L = 1$, the term $II$ is missing.

Consider $I$. Let $\frac{1}{\Theta} \geq \frac{1}{r}$ and $0 < t \leq 1$. Then $t^{\frac{1}{\Theta}-\frac{1}{r}} \leq 1$. Hence, using the properties (P2) and (P4), we immediately get

$$I \leq \rho_E^{\kappa_X}\left(\chi_{(0,1]}(t)\right) < \infty.$$

Note that here $q_X$ may also take the value $\infty$. The case $q_X < \infty$ and $\frac{1}{\Theta}+\frac{1}{q_X} > \frac{1}{r}$ is considered in [12, proof of Theorem 4.1]. (Note that these two cases overlap.) So, it has been proven that $X_0 \cap X_1 \hookrightarrow \bar{X}_{\Theta,r;X}$ and $(X_0, X_1)_{\Theta,r;X}$ is an intermediate space for the couple $\bar{X}$.

Step 3. Here we will show that $\bar{X}_{\Theta,r;E}$ is an exact interpolation functor. Suppose that the linear operator $T: X_0 + X_1 \to X_0 + X_1$ is such that

$$\|Tx\|_{X_i} \le M_i \|x\|_{X_i}, \qquad i = 0,1.$$

By the definition of the *K*-functional, we obtain that

$$K\left(t^{\frac{1}{\Theta}}, T(x)\right) \le \inf\left\{\|T(x_0)\|_{X_0} + t^{\frac{1}{\Theta}}\|T(x_1)\|_{X_1} : x = x_0 + x_1, x_i \in X_i, i = 0,1\right\}$$

$$\le \inf\left\{M_0\|x_0\|_{X_0} + t^{\frac{1}{\Theta}} M_1\|x_1\|_{X_1} : x = x_0 + x_1, x_i \in X_i, i = 0,1\right\}$$

$$\le \max\{M_0, M_1\} \inf\left\{\|x_0\|_{X_0} + t^{\frac{1}{\Theta}}\|x_1\|_{X_1} : x = x_0 + x_1, x_i \in X_i, i = 0,1\right\}$$

$$= \max\{M_0, M_1\} K\left(t^{\frac{1}{\Theta}}, x\right).$$

Hence,

$$\|Tx\|_{\bar{X}_{\Theta,r;X}} = \rho_E\left(t^{-\frac{1}{r}} K\left(t^{\frac{1}{\Theta}}, T(x)\right)\right)$$

$$\le \max\{M_0, M_1\}\, \rho_E\left(t^{-\frac{1}{r}} K\left(t^{\frac{1}{\Theta}}, x\right)\right) = \max\{M_0, M_1\}\, \|x\|_{\bar{X}_{\Theta,r;X}}.$$

So, $(*,*)_{\Theta,r;X}$ is an exact interpolation functor. ■

## 4. Main result

Now our main theorem can be proved.

Theorem 15. Let $E$ be a r.i. quasi-Banach function space on $(0, L)$, $0 < p_0 < p_1 \le \infty$ and $0 < q_0, q_1 \le \infty$. Assume that $p_0 < p_E$ and

(i) either $q_E < p_1 < \infty$ and $\frac{1}{\Theta} = \frac{1}{p_0} - \frac{1}{p_1}$

(ii) or $p_1 = q_1 = \infty$, $q_E \le \infty$ and $\Theta = p_0$.

Then

$$\left(L_{p_0,q_0}(0, L), L_{p_1,q_1}(0, L)\right)_{\Theta,p_0,E} = E.$$

Proof. Put $F = \left(L_{p_0,q_0}, L_{p_1,q_1}\right)_{\Theta,p_0,E}$. We must prove that $F = E$. As usual, $\rho_E$ denotes the quasi-norm on $E$. Let $f \in L_{p_0,q_0} + L_{p_1,q_1}$.

Case $q_E < p_1 < \infty$ and $\frac{1}{\Theta} = \frac{1}{p_0} - \frac{1}{p_1}$. Holmstedt's formula [17, Theorem 4.2] gives us

$$K\left(t^{\frac{1}{\Theta}}, f; L_{p_0,q_0}, L_{p_1,q_1}\right) \sim \left(\int_0^t \left\{s^{\frac{1}{p_0}} f^*(s)\right\}^{q_0} \frac{ds}{s}\right)^{\frac{1}{q_0}} + t^{\frac{1}{\Theta}}\left(\int_t^L \left\{s^{\frac{1}{p_1}} f^*(s)\right\}^{q_1} \frac{ds}{s}\right)^{\frac{1}{q_1}}, \quad t \in (0, L).$$

Thus,

$$t^{-\frac{1}{p_0}} K\left(t^{\frac{1}{\Theta}}, f; L_{p_0,q_0}, L_{p_1,q_1}\right)$$

$$\sim t^{-\frac{1}{p_0}}\left(\int_0^t \left\{s^{\frac{1}{p_0}} f^*(s)\right\}^{q_0} \frac{ds}{s}\right)^{\frac{1}{q_0}} + t^{-\frac{1}{p_1}}\left(\int_t^L \left\{s^{\frac{1}{p_1}} f^*(s)\right\}^{q_1} \frac{ds}{s}\right)^{\frac{1}{q_1}}$$

$$= \left(H^{(p_0,q_0)} f\right)(t) + \left(H_{(p_1,q_1)} f\right)(t)$$

and

$$\|f\|_F = \rho_E\left(t^{-\frac{1}{p_0}} K\left(t^{\frac{1}{\Theta}}, f; L_{p_0,q_0}, L_{p_1,q_1}\right)\right)$$

$$\sim\rho_E\left(\left(H^{(p_0,q_0)}f\right)(t)\right)+\rho_E\left(\left(H_{(p_1,q_1)}f\right)(t)\right).$$

Because $p_0 < p_E$ and $q_E < p_1 < \infty$, according to Theorem 11, we conclude

$$\|f\|_F\sim\|f\|_E.$$

Case $q_E \leq \infty$, $p_1 = q_1 = \infty$, $\Theta = p_0$. In this case, the proof is simpler because

$$t^{-\frac{1}{p_0}}K\left(t^{\frac{1}{p_0}},f;L_{p_0,q_0},L_\infty\right)\sim t^{-\frac{1}{p_0}}\left(\int_0^t\left\{s^{\frac{1}{p_0}}f^*(s)\right\}^{q_0}\frac{ds}{s}\right)^{\frac{1}{q_0}}=\left(H^{(p_0,q_0)}f\right)(t),\quad t\in(0,L).$$

Hence,

$$\|f\|_F\sim\rho_E\left(\left(H^{(p_0,q_0)}f\right)(t)\right)\sim\|f\|_E. \qquad \blacksquare$$

## 5. *Important corollaries*

In this section, some corollaries of Theorem 15 for the case $q_0 = p_0$ and $q_1 = p_1$ will be formulated and compared with some known results. The following corollary coincides with Corollary 4.3 in [12] if $q_E < p_1 < \infty$.

Corollary 16. Let $E$ be a r.i. quasi-Banach function space on $(0,L)$ and $0 < p_0 < p_1 \leq \infty$. Assume that $p_0 < p_E$ and

(i) either $q_E < p_1 < \infty$ and $\frac{1}{\Theta} = \frac{1}{p_0} - \frac{1}{p_1}$

(ii) or $p_1 = \infty$, $q_E \leq \infty$ and $\Theta = p_0$.

Then

$$\left(L_{p_0}(0,L),L_{p_1}(0,L)\right)_{\Theta,p_0,E}=E. \tag{2}$$

In particular, (Cf. [3, Theorem II.6.6 and Theorem III.2.2] and [21, Theorem II.4.1].)

$$L_{p_0}(0,L)\cap L_{p_1}(0,L)\hookrightarrow E\hookrightarrow L_{p_0}(0,L)+L_{p_1}(0,L).$$

Remark 17. Theorem 1 in the pioneering work of D. Boyd [5] gives sufficient conditions for a r.i. Banach function space to be an interpolation space between $L_{p_0}$ and $L_{p_1}$ if $0 < p_0 < p_1 < \infty$. Many papers are devoted to this problem. We will only mention [18, Theorem 73] and [7, Propositions 2.13 and 2.15]. Note that from (2) and Theorem 14, it follows that Corollary 16 improves and extends Theorem 1 from [5].

The following example shows that the condition $p_0 < p_E$ in Corollary 16 cannot be improved.

Example 18. Let us take any $q \in (0,1)$ and consider the space $E = L_{1,q}(0,\infty)$. It is not hard to check that

$$p_E = q_E = 1.$$

It follows from Corollary 16 that for each $p_0 \in (0,1)$

$$L_{p_0}(0,\infty)\cap L_\infty(0,\infty)\hookrightarrow E.$$

We will show that

$$L_1(0,\infty)\cap L_\infty(0,\infty)\not\subseteq E.$$

We will use Lorentz sequence space $l_{1,q}$. This is the set of all complex-valued sequences $\{f_k\}$ $(k = 1,2,\ldots)$ such that

$$\|\{f_k\}\|_{1,q}:=\left\{\sum_{k=1}^\infty\left(k^{1-\frac{1}{q}}f_k^*\right)^q\right\}^{\frac{1}{q}}<\infty,$$

where $\{f_k^*\}$ $(k = 1,2,\ldots)$ is the non-increasing rearrangement of the sequence $\{|f_k|\}$. It is known [2] that

$$l_{1,q} \subsetneq l_1.$$

Hence, there exists a positive, non-increasing sequence $\{f_k\}$ $(k = 1,2,\dots)$ such that

$$\|\{f_k\}\|_1 = \sum_{k=1}^{\infty} f_k < \infty \qquad \text{and} \qquad \|\{f_k\}\|_{1,q} = \left\{\sum_{k=1}^{\infty}\left(k^{1-\frac{1}{q}} f_k\right)^q\right\}^{\frac{1}{q}} = \infty.$$

Consider the function

$$f(t) := f_k \qquad \text{if } k-1 < t \le k, \qquad (k = 1,2,\dots).$$

Obviously, $f^* = f$, $f \in L_\infty$ and

$$\|f\|_1 = \sum_{k=1}^{\infty} \int_{k-1}^{k} f(t)dt = \sum_{k=1}^{\infty} f_k < \infty.$$

Thus,

$$f \in L_1 \cap L_\infty.$$

Furthermore, since $t^{q-1}$ decreases

$$\left(\|f\|_{1,q}\right)^q = \sum_{k=1}^{\infty} \int_{k-1}^{k} \left(t^{1-\frac{1}{q}} f(t)\right)^q dt = \sum_{k=1}^{\infty} (f_k)^q \int_{k-1}^{k} t^{q-1} dt$$

$$\ge \sum_{k=1}^{\infty} (f_k)^q k^{q-1} = \sum_{k=1}^{\infty}\left(k^{1-\frac{1}{q}} f_k\right)^q = \infty.$$

Thus,

$$f \notin L_{1,q} = E.$$

Corollary 19. Let $E$ be an arbitrary r.i. quasi-Banach function space on $(0, L)$. Then there exist $p_0$ and $p_1$ such that $0 < p_0 < p_1 \le \infty$ and

$$\left(L_{p_0}(0,L), L_{p_1}(0,L)\right)_{\Theta,p_0,E} = E, \qquad \text{where} \qquad \frac{1}{\Theta} = \frac{1}{p_0} - \frac{1}{p_1}.$$

In addition, there exists $p_0$ such that $0 < p_0 < \infty$ and

$$\left(L_{p_0}(0,L), L_\infty(0,L)\right)_{p_0,p_0,E} = E.$$

Proof. If $q_E < \infty$, one must take any $p_0, p_1$ such that $0 < p_0 < p_E$ and $q_E < p_1 \le \infty$. If $q_E = \infty$, one must take any $p_0$ such that $0 < p_0 < p_E$ and $p_1 = \infty$. ■

Corollary 16 has the following equivalent formulation.

Corollary 20.

(i) Let $0 < p_0 < p_1 < \infty$. Then for each r.i. quasi-Banach function space $E$ on $(0, L)$ such that $p_0 < p_E \le q_E < p_1$ holds

$$\left(L_{p_0}(0,L), L_{p_1}(0,L)\right)_{\Theta,p_0,E} = E,$$

where $\frac{1}{\Theta} = \frac{1}{p_0} - \frac{1}{p_1}$.

(ii) Let $0 < p_0 < \infty$. Then for each r.i. quasi-Banach function space $E$ on $(0, L)$ such that $p_0 < p_E$ holds

$$\left(L_{p_0}(0,L), L_\infty(0,L)\right)_{p_0,p_0,E} = E.$$

The following statement is probably well known. But it is quite difficult to find a reference for it.

Corollary 21. The set of all countable-valued functions[2] is dense in any r.i. quasi-Banach function space $E$.

---

[2] A countable-valued function is a function whose range is a countable set.

Proof. By Corollary 19, there exists $0 < p_0 < \infty$ such that $L_{p_0} \cap L_\infty \hookrightarrow E \hookrightarrow L_{p_0} + L_\infty$. The further proof is as in the proof of Lemma II.4.2 in [21]. ∎

## Acknowledgements

The author would like to thank Prof. Lech Maligranda for pointing out to us the papers [19, 20, 22 and 26].